\documentclass[11pt, english]{article}
\usepackage[cp866]{inputenc}
\usepackage{babel}
\usepackage{amsfonts}
\usepackage{amssymb}
\usepackage{amsthm}
\usepackage{amsmath}
\usepackage{enumerate}
\usepackage{mathrsfs}

\newtheorem{thm}{Theorem}
\newtheorem{prop}{Proposition}
\newtheorem{cor}{Corollary}

\newtheorem{defn}{Definition}
\newtheorem{fact}{Fact}
\newtheorem{claim}{Claim}
\newtheorem{quest}{Question}

\newcommand{\A}{\mbox{$\mathcal{A}$}}

\newcommand{\LL}{\mbox{$\mathcal{L}$}}

\newcommand{\at}{\char'100}

    \DeclareMathOperator{\lh}{lh}

    \DeclareMathOperator{\dom}{dom}
    
    \DeclareMathOperator{\ran}{ran}

     \DeclareMathOperator{\II}{II}

    \DeclareMathOperator{\I}{I}

    \DeclareMathOperator{\low}{low}
    \DeclareMathOperator{\ck}{CK}
    \DeclareMathOperator{\codom}{codom}
    \DeclareMathOperator{\diam}{diam}

    \newcommand{\restrict}{\upharpoonright}

    \newcommand{\baire}{\mathscr N}
    \newcommand{\forces}{\Vdash}

\def\bSigma{{\boldsymbol{\Sigma}}}
\def\bPi{{\boldsymbol{\Pi}}}

\title{The Effective Theory of Borel Equivalence Relations\thanks{The
authors acknowledge the generous support of the FWF through project
P 19375-N18}}
\author{Ekaterina B. Fokina\and Sy-David Friedman\and Asger T\"{o}rnquist}

\begin{document}
\selectlanguage{english}
\maketitle

\begin{abstract}
The study of Borel equivalence relations under Borel reducibility
has developed into an important area of descriptive set theory. The
dichotomies of Silver (\cite{silver80}) and
Harrington-Kechris-Louveau (\cite{hakelou90}) show that with respect
to Borel reducibility, any Borel equivalence relation strictly above
equality on $\omega$ is above equality on ${\cal P}(\omega)$, the
power set of $\omega$, and any Borel equivalence relation strictly
above equality on the reals is above equality modulo finite on
${\cal P}(\omega)$. In this article we examine the effective content
of these and related results by studying effectively Borel
equivalence relations under effectively Borel reducibility. The
resulting structure is complex, even for equivalence relations with
finitely many equivalence classes. However use of Kleene's $O$ as a
parameter is sufficient to restore the picture from the noneffective
setting. A key lemma is the existence of two effectively Borel sets
of reals, neither of which contains the range of the other under any
effectively Borel function; the proof of this result applies Barwise
compactness to a deep theorem of Harrington (see \cite{hjohar})
establishing for any recursive ordinal $\alpha$ the existence of
$\Pi^0_1$ singletons whose $\alpha$-jumps are Turing incomparable.
\end{abstract}

\section{Introduction}

If $E$ and $F$ are Borel equivalence relations on Polish spaces $X$
and $Y$ respectively, then $E$ is \emph{Borel reducible} to $F$ if
and only if there is a Borel function $f:X\to Y$ such that $xEy$ if
and only if $f(x)Ff(y)$. The study of Borel equivalence relations
under Borel reducibility has developed into a rich area of
descriptive set theory. Surveys of some of this work may be found in
\cite{beke96, gao09, hjorth, hjoke01, kanovei08, Kechris99}. In the
noneffective setting, Borel equivalence relations with countably
many equivalence classes are equivalent (i.e. bi-reducible) exactly
if they have the same number of equivalence classes. For Borel
equivalence relations with uncountably many equivalence classes
there are two fundamental dichotomies:

\medskip

\noindent \textbf{The Silver Dichotomy} (\cite{silver80}). \emph{If
$E$ is a Borel equivalence relation with uncountably many
equivalence classes then equality on ${\cal P}(\omega)$, the power
set of $\omega$, is Borel reducible to $E$}.

\medskip

\noindent \textbf{The Harrington-Kechris-Louveau Dichotomy}
(\cite{hakelou90}). \emph{If $E$ is a Borel equivalence relation not
Borel reducible to equality on ${\mathcal P}(\omega)$ then $E_0$ is
Borel reducible to $E$, where $E_0$ is equality modulo finite on
${\cal P}(\omega)$.}

\medskip

In this article we introduce the effective version of this theory.
If $E$ and $F$ are effectively Borel (i.e., $\Delta^1_1$)
equivalence relations on effectively presented Polish
spaces\footnote{In the sense of Moschovakis,
\cite[3B]{moschovakis80}. In this paper we will deal almost
exclusively with the spaces $\omega$, $\mathcal P(\omega)$ and
$\baire=\omega^\omega$.} spaces $X$ and $Y$, respectively, then we
say that $E$ is \emph{effectively Borel reducible} to $F$ if there
is an effectively Borel function $f:X\to Y$ such that $xEy$ if and
only if $f(x)Ff(y)$. The resulting effective theory reveals an
unexpectedly rich new structure, even for equivalence relations with
finitely many classes. For $n\leq\omega$, let $=_n$ denote equality
on $n$, let $=_{{\mathcal P}(\omega)}$ denote equality on the power
set of $\omega$ and let $E_0$ denote equality modulo finite on
${\mathcal P}(\omega)$. The notion of effectively Borel reducibility
on effectively Borel equivalence relations naturally gives rise to a
degree structure, which we denote by ${\cal H}$.

We show the following:

\medskip

\noindent \textbf{Theorem A.} \emph{For any finite $n$, the partial order of $\Delta^1_1$ subsets of $\omega$ under inclusion can be order-preservingly embedded into $\cal H$ between the degrees of $=_n$ and $=_{n+1}$.
The same holds between the degrees of $=_\omega$ and $=_{{\mathcal P}(\omega)}$, and between
$=_{{\mathcal P}(\omega)}$ and $E_0$.}

\medskip

A basic tool in the proof of Theorem A is the following result,
which may be of independent interest:

\medskip

$(*)$ \emph{There are effectively Borel sets $A$ and $B$ such that
for no effectively Borel function $f$ does one have $f[A]\subseteq
B$ or $f[B]\subseteq A$.}

\medskip

$(*)$ is proved via a Barwise compactness argument applied to a deep
result of Harrington (see \cite{hjohar}) establishing for any
recursive ordinal $\alpha$ the existence of $\Pi^0_1$ singletons
whose $\alpha$-jumps are Turing incomparable.

We also examine the effectivity of the Silver and
Harrington-Kechris-Louveau dichotomies. Harrington's proof of the
Silver dichotomy (see \cite{gao09} or \cite{jech03}) and the
original proof of the Harrington-Kechris-Louveau dichotomy in
\cite{hakelou90} respectively show that if an effectively Borel
equivalence relation has countably many equivalence classes then it
is effectively Borel reducible to $=_\omega$ and if it is Borel
reducible to $=_{{\cal P}(\omega)}$ then it is in fact effectively
Borel reducible to $=_{{\cal P}(\omega)}$. We complete the picture
by showing:

\medskip

\noindent \textbf{Theorem B.} \emph{Let $O$ denote Kleene's $O$. If
an effectively Borel equivalence relation $E$ has uncountable many
equivalence classes then there is a $\Delta^1_1(O)$ function
reducing $=_{{\cal P}(\omega)}$ to $E$, and this parameter is best
possible. If an effectively Borel equivalence relation $E$ is not
Borel reducible to $=_{{\mathcal P}(\omega)}$ then there is a
$\Delta^1_1(O)$ function reducing $E_0$ to $E$, and this parameter
is best possible.}

\medskip

In other words, while Theorem A rules out that the dichotomy
Theorems of Silver and Harrington-Kechris-Louveau are effective,
Theorem B shows that the Borel reductions obtained in the dichotomy
Theorems can in fact be witnessed by $\Delta^1_1(O)$ functions, and
that Kleene's $O$ is the best possible parameter we can hope for in
general. The proof of Theorem B is based on a detailed analysis of
the effectiveness of category notions in the Gandy-Harrington
topology, due to the third author.

There remain many open questions in the effective theory. We mention
a few of them at the end of the article.

\medskip

\noindent{\bf Organization.} The paper is organized into 6 sections.
In \S 2 we introduce some basic notation used in the paper, and
recall some well-known theorems and facts that our proofs rely on.
In \S 3 we prove $(*)$, which serves as a basic tool throughout the
paper. The proof of Theorem A and several extensions of Theorem A is
found in \S 4. In \S 5 we give a detailed analysis of the
effectiveness of category notions in the Gandy-Harrington topology.
Finally, Theorem B is proved in \S6.

\section{Background}

Throughout this paper, Hyp stands for $\Delta^1_1$, both for subsets
of $\omega$ and for subsets of Baire space $\baire=\omega^{\omega}$.
Elements of $\baire$ are called ``reals''. We state without proofs
some well-known results that we will need in this paper. For further
details the reader may consult the provided references.

For a linear ordering $<$ denote by $\mathcal{W}f(<)$ the largest
well-ordered initial segment of $<$. We can identify
$\mathcal{W}f(<)$ with an ordinal without danger of confusion.

\begin{thm}[Barwise, see \cite{Barwise}]\label{Barwise}
Let $\LL$ be a recursive language, $\A=L_{\omega_1^{\ck}}$, and let
$\LL_{\mathcal{A}}$ be $\LL_{\omega_1\omega}$ restricted to
$\varphi\in \A$. Suppose $\Phi\subseteq\LL_{\mathcal{A}}$ is a
$\Sigma_1(\A)$ set of sentences and every $\Phi_0\subseteq\Phi$ such
that $\Phi_0\in \A$ has a model. Then $\Phi$ has a model. Moreover,
if $<\in\LL$ and for all $\alpha<\omega_1^{\ck}$ and
$\Phi_0\subseteq\Phi$ such that $\Phi_0\in \A$ there is a model of
$\Phi_0$ in which $<$ is a linear ordering of length at least
$\alpha$, then $\Phi$ has a model in which $<$ is a linear ordering
satisfying $\mathcal{W}f(<)=\omega_1^{\ck}\neq Field(<)$.
\end{thm}

\begin{defn}
Let $\Gamma$ be a point-class (in the sense of Moschovakis
\cite{moschovakis80}) and let $A$ be a set of reals. We call $A$
\emph{a $\Gamma$ singleton} iff $A$ has exactly one element and $A$
belongs to $\Gamma$.
\end{defn}

In this paper $\Gamma$ will usually be $\Pi^0_1$ or $\Delta^1_1$ (i.e. Hyp).

\begin{fact}[see \cite{Rogers,Sacks}]\label{fct0}
 \begin{enumerate}
  \item Every Hyp real is a Hyp singleton;
  \item A countable Hyp set of reals contains only Hyp reals;
  \item For every Hyp real $X$ there is a $\Pi^0_1$ singleton $Y$, such that $X\leq_T Y$.
 \end{enumerate}
\end{fact}

\begin{thm}[Uniformization, see e.g. Chapter II of \cite{Sacks}]\label{Uni}
 Every $\Pi^1_1$ binary relation on $\baire\times \baire$ contains a $\Pi^1_1$
 function with the same domain.
\end{thm}

\begin{thm}[Dependent Choice, see Chapter II of \cite{Sacks}]\label{DC}
If $P$ is a Hyp binary relation and for all Hyp reals $X$ there
exists a Hyp real $Y$ such that $P(X,Y)$, then for all Hyp reals
$X$, there is a Hyp $\omega$-sequence $X=X_0, X_1,\ldots$ such that
$P(X_n,X_{n+1})$, for all $n$.
\end{thm}

Recall that $E_0$ is the equivalence relation on $2^\omega$
defined by
$$
xE_0 y\iff (\exists n)(\forall m\geq n) x(m)=y(m),
$$
equivalently, $E_0$ is equality modulo finite in $\mathcal
P(\omega)$. The next result is folklore (see e.g. \cite[Theorem
3.2]{hjorth}):

\begin{fact}\label{fct1} If $h: 2^\omega\to 2^\omega$ is Baire measurable and constant on $E_0$
classes then $h$ is constant on a comeagre set.
\end{fact}

The following result will be used several times:

\begin{fact}[Kechris, \cite{kechris73}]\label{fct2} If $B\subset 2^\omega\times 2^\omega$ is Hyp then $\{x: \{y: (x,y)\in
B\} \text{ is non-meagre}\}$ is $\Sigma^1_1$.
\end{fact}

Finally, we will use the following result from \cite{Harrington}.
For a sketch of the proof see also \cite{hjohar}.

\begin{thm}\label{Har}
For any recursive ordinal $\alpha$ there is a sequence of reals
$\langle a_n | n < \omega\rangle$ such that for some recursive
sequence $\langle\varphi_n | n < \omega\rangle$ of $\Pi^0_1$
formulas, $a_n$ is the unique solution to $\varphi_n$ for each $n$
and no $a_n$ is recursive in the $\alpha$-jump of $\langle a_m | m
\neq n\rangle$.
\end{thm}

\noindent \textbf{Remark.} We will also use the following weaker
form of Theorem \ref{Har}. For every recursive ordinal $\alpha$
there are two $\Pi^0_1$ singletons $a, b$ such that $a$ is not
recursive in the $\alpha$-jump of $b$ and $b$ is not recursive in
the $\alpha$-jump of $a$.

\medskip

\noindent \textbf{Notation.} If $a$ is a real and
$\alpha<\omega_1^{\ck}$ then we denote by $a^{\alpha}$ the
$\alpha$-jump of $a$.

\section{The basic tool: Hyp incomparable Hyp sets of reals }

The theorem which we prove in this section will be used repeatedly
to obtain the results of this paper.

\begin{thm}\label{MainThm}
There exist two nonempty $\Pi^0_1$ sets $A,B\subseteq\baire$, such
that for no Hyp function $F:\baire\to\baire$ do we have
$F[A]\subseteq B$ or $F[B]\subseteq A$.
\end{thm}

\noindent \textbf{Remark.} If $A$ and $B$ are as in Theorem
\ref{MainThm} then neither $A$ nor $B$ contains a Hyp real: Suppose
$A$ contains a Hyp real $y$; then the constant function with value
$y$ maps $B$ into $A$, contradiction. In particular, it follows that
there is no Hyp $F$ such that $F[\sim A]\subseteq B$ or $F[\sim
B]\subseteq A$.

\begin{proof}[Proof of Theorem \ref{MainThm}]
Let $\A=L_{\omega_1^{\ck}}$,
$\LL\supseteq\{\in,<,\underline{x_0},\underline{x_1}\}\cup\{\underline{\alpha}:\alpha\in
\A\}$, where $\underline{x_0},\underline{x_1}$ and
$\underline{\alpha}$ are constant symbols. Consider the set of
sentences $\Phi$ consisting of:
\begin{enumerate}[(1)]
 \item $ZF^-$
 \item ($\forall x)(x\in \underline{\omega}\leftrightarrow \bigwedge_n x=\underline{n})$
 \item $<=\in\upharpoonright\text{Ordinals}$
 \item $\underline{x_0},\underline{x_1}\subseteq\underline\omega$
 \item $\bigvee_{\varphi\in\Pi^0_1}\left[(\exists!v)\varphi(v)\wedge\varphi(\underline{x_i})\right]$\ ($i=0,1$, $\varphi$ ranges over all $\Pi^0_1$
 formulas.)
 \item $\underline{x_0}\nleq_T \underline{x_1}^{\underline{\alpha}}, \underline{x_1}\nleq_T \underline{x_0}^{\underline{\alpha}}\text{, for all }\alpha\in \omega_1^{\ck}\text{.}$

\end{enumerate}
The set $\Phi$ is a $\Sigma_1$ set of sentences. By the remark
following Theorem \ref{Har}, for every recursive ordinal $\alpha$
there exist $\Pi^0_1$ singletons $a_{\alpha},b_{\alpha}$, such that
$a_{\alpha}$ is not recursive in the $\alpha$th Turing jump of
$b_{\alpha}$ and $b_{\alpha}$ is not recursive in the $\alpha$th
Turing jump of $a_{\alpha}$. Thus, we can apply Theorem
\ref{Barwise}. We get a model $\langle
M,E,<,x_0,x_1\rangle\models\Phi$ such that
$L_{\omega_1^{\ck}}\subseteq M$, $M$ has nonstandard ordinals and
every standard ordinal of $M$ is recursive, i.e., the standard part
of $<^M$ is $\omega_1^{\ck}$. Then in $M$ there must be $\Pi^0_1$ singletons $a$ and $b$
such that $a\nleq_T b^{\alpha}, b\nleq_T
a^{\alpha}$ for $\alpha<\omega_1^{\ck}$ and since
$\omega_1^a=\omega_1^b=\omega_1^{\ck}$, $a$ and $b$ are
Hyp-incomparable.

Choose $\Pi^0_1$ formulas $\varphi_a$ and $\varphi_b$, such that in
$M$, $\varphi_a(x)\leftrightarrow x=a$ and
$\varphi_b(x)\leftrightarrow x=b$. Note that $a$ and $b$ are the
unique solutions of $\varphi_a$ and $\varphi_b$ in $M$,
respectively. Then the formulas $\varphi_a$ and $\varphi_b$ define
$\Pi^0_1$ sets (not singletons) in ${V}$. Let $A=\{x:\varphi_a(x)\}$
and $B=\{x:\varphi_b(x)\}$.

\begin{claim}\label{IncSets}
There is no Hyp function $F$ such that $F[A]\subseteq B$.
\end{claim}

\begin{proof}
Suppose $F$ were such a function. Consider $F(a)\in M$. It is Hyp in
$a$. On the other hand, $F(a)=t\in B$. Therefore by definition of
$B$, $\varphi_b(t)$ holds in $M$, and so $t=b$. Thus, $b$ is Hyp in
$a$, contradicting the properties of $a$ and $b$.
\end{proof}

This completes the proof of the theorem.
\end{proof}

\begin{thm}\label{MainInf}
There exists a uniform sequence $A_0, A_1,\ldots$ of nonempty
$\Pi^0_1$ sets such that for each $n$ there is no Hyp function $F$
such that $F[A_n]\subseteq \bigcup_{m\neq n}A_m$.
\end{thm}

\begin{proof}
The proof is analogous to the previous proof using Theorems
\ref{Har} and \ref{Barwise}. We consider
$\LL\supseteq\{\in,<,\underline{x_0},\underline{x_1},\ldots\}\cup\{\underline{\alpha}:\alpha\in
\A\}$ and the following set of sentences $\Phi$:
\begin{enumerate}[(1)]
 \item $ZF^-$
 \item $(\forall x)(x\in \underline{\omega}\leftrightarrow \bigwedge_n x=\underline{n})$
 \item $<=\in\upharpoonright\text{Ordinals}$
 \item $\bigwedge_n\underline{x_n}\subseteq\underline\omega$
 \item $\bigwedge_n\bigvee_{\varphi\in\Pi^0_1}\left[(\exists!v)\varphi(v)\wedge\varphi(\underline{x_n})\right]$,\ ($\varphi$ ranges over all $\Pi^0_1$
 formulas)
 \item $\bigwedge_{m\neq n} \underline{x_m}\nleq_T \underline{x_n}^{\underline{\alpha}}\text{ for all }\alpha\in \omega_1^{\ck}\text{.}$

\end{enumerate}
By the properties of sequences $\langle a_n : n < \omega\rangle$
from Theorem \ref{Har}, we get that the resulting sequence
$A_0,A_1,\ldots$ of $\Pi^0_1$ sets is uniform and has the required
properties exactly as in Theorem \ref{MainThm}.
\end{proof}

\section{Hyp Equivalence Relations under Hyp Reducibility}

\begin{defn}
Let $E$ and $F$ be equivalence relations on $\baire$. We say that
$E$ is \emph{Hyp-reducible} to $F$ if there exists a Hyp function
$f:\baire\to\baire$ such that
$$
xEy\iff f(x)Ff(y),
$$
in which case we will write $E\leq_H F$.
\end{defn}
This notion induces a natural notion of \emph{Hyp-equivalence} (or
\emph{Hyp bi-reducibility}) and \emph{Hyp-degrees}: we let
$E\equiv_H F$ if and only if $E\leq_H F$ and $F\leq_H E$.

\begin{defn}
For every $n\in\omega, n\geq 1$, let $=_n$ be the Hyp-degree of the following
equivalence relation:
$$
x\equiv y\Longleftrightarrow x(0)=y(0) \text{ or both }
x(0),y(0)\geq n-1.
$$
The Hyp-degree {$=_{\omega}$} is the Hyp-degree of the equivalence
relation
$$
x\equiv y\Longleftrightarrow x(0)=y(0).
$$
\end{defn}

\subsection{Hyp Equivalence Relations with countably many classes}

\begin{prop}\label{prop1}
Let $1\leq n\leq\omega$ and let $E$ be a Hyp equivalence relation. Then
$=_n\leq_HE$ iff $E$ has at least $n$ classes containing Hyp reals.
\end{prop}

\begin{proof}
$(\Rightarrow):$ For every $1\leq n\leq\omega$, the equivalence relation
$=_n$ has exactly $n$ equivalence classes and each of them contains
a Hyp real. Under Hyp-reducibility Hyp reals are sent to Hyp reals,
equivalent reals are sent to equivalent reals, non-equivalent reals
are sent to non-equivalent reals.

$(\Leftarrow):$ If $n$ is finite, pick $n$ Hyp reals
$x_0,\ldots,x_{n-1}$ that lie in different equivalence classes of
$E$. The function $F$ that sends the $i$th equivalence class of
$=_n$ to $x_i$ witnesses the reduction. To prove the result for
$n=\omega$ we use Theorem \ref{DC}. Suppose $E$ is an equivalence
relations with infinitely many classes containing Hyp reals. We want
to prove that $=_{\omega}$ Hyp-reduces to $E$. We will find a Hyp
sequence of equivalence classes of $E$ with Hyp reals in them.
Consider the following relation $P(X,Y)$ on
$\omega\times\baire^{<\omega}$:
\begin{align*}
P(X,Y)\iff & [X=(n,X_0,\ldots,X_n)\wedge\bigwedge_{i\neq j} \neg
X_iEX_j]\longrightarrow\\
& [Y=(n+1,Y_1,\ldots,Y_n,Y_{n+1})\wedge\bigwedge_i
X_i=Y_i\wedge\bigwedge_{i\neq j}\neg Y_iEY_j]
\end{align*}

Then $P$ is Hyp. Moreover, as $E$ has infinitely many Hyp classes,
for every Hyp $X$ there exists a Hyp $Y$ such that $P(X,Y)$. It
follows from Theorem \ref{DC} that there exists a uniform sequence
of Hyp sets $X_0,X_1,\ldots$ such that
$$
\forall i,j(i\neq j\rightarrow \neg X_iEX_j).
$$
Then the function that sends the equivalence class $\{x:x(0)=n\}$ of
$=_{\omega}$ to $X_n$ is Hyp and witnesses the reduction.
\end{proof}

\begin{cor}
If $=_n\leq_H E$, for all $1\leq n<\omega$, then $=_{\omega}\leq_H E$.
\end{cor}

\begin{prop}\label{prop2}
Let $1\leq n\leq\omega$ and let $E$ be a Hyp equivalence relation. Then
$E\leq_H =_n$ iff $E$ has at most $n$ classes.
\end{prop}

\begin{proof}
The direction $(\Rightarrow)$ is obvious since non-equivalent reals
are sent to non-equivalent reals under Hyp-reducibility.

To prove $(\Leftarrow)$ we need to show that the equivalence classes
of a Hyp equivalence relation with at most countably many
equivalence classes are uniformly Hyp.


By Harrington's proof of the Silver Dichotomy (see \cite[Theorem
32.1]{jech03} or \cite[Theorem 5.3.5]{gao09}), if $E$ has only countably many classes then every real belongs
to a Hyp subset of some equivalence class. Let $C$ be the set of
codes for Hyp subsets of an equivalence class; then $C$ is
$\Pi^1_1$. Consider the relation
$$
R = \{(x,c) : c\in C \text{ and } x\in H(c)\text{, the Hyp set coded
by }c\}.
$$
Then $R$ is $\Pi^1_1$ and can be uniformised by a $\Pi^1_1$ function
$F$. As the values of $F$ are numbers, $F$ is Hyp and by separation
we can choose a Hyp $D\subseteq C$, $D\supseteq \ran(F)$. Now define
an equivalence relation $E^*$ on $D$ by:
\begin{align*}
 d_0 E^* d_1 & \iff (\forall x_0,x_1) (x_0\in H(d_0)\wedge x_1\in H(d_1))\rightarrow x_0Ex_1\\
 & \iff (\exists x_0,x_1) (x_0\in H(d_0)\wedge x_1\in H(d_1)\wedge x_0Ex_1).
\end{align*}
i.e. $d_0 E^* d_1$ if and only if $H(d_0)$  and $H(d_1)$ are
subsets of the same $E$-equivalence class. Note that $E^*$ is Hyp.
The relation $E$ Hyp-reduces to $E^*$ via $x\mapsto F(x)$. But $E^*$
is just a Hyp relation on a Hyp set of numbers, so $E^*$ is
Hyp-reducible to $=_{\omega}$ (to see this, send $c$ to the least
number $c^*$, $cE^*c^*$).

Thus if $E$ is a Hyp equivalence relation with at most countably
many classes then $E$ is Hyp-reducible to $=_{\omega}$. (In
particular, all equivalence classes of $E$ are Hyp.) One can
similarly see that if $E$ has at most $n$ classes then $E$ is
Hyp-reducible to $=_n$.
\end{proof}


Obviously, the degree $=_1$ is Hyp-reducible to any other
Hyp-degree. But $=_2$, the equivalence relation with the 2 classes
$\{x:x(0)=0\}$ and $\{x:x(0)\geq 1\}$ is \emph{not} the successor to
$=_1$. This is the content of the next theorem.

\begin{thm}\label{Thm1-2}
\begin{enumerate}
 \item There is a Hyp equivalence relation strictly between $=_1$ and~$=_2$.
 \item For every finite $n$, there is a Hyp equivalence relation strictly between $=_n$ and $=_{n+1}$.
 \item For every $n_0<n_1\leq\omega$, there is a Hyp equivalence relation above $=_{n_0}$, below $=_{n_1}$ and incomparable with $=_n$, for all $n_0<n<n_1$.
\end{enumerate}

\end{thm}

\begin{proof} The proof is based on the following fact.
\begin{fact}[\cite{Sacks}]\label{HNH}
There is a nonempty Hyp set $X$ which contains no Hyp reals.
\end{fact}

To prove the first statement, take a Hyp equivalence relation $E$
with two equivalence classes $X$ and $\sim X$, where $X$ is from
Fact \ref{HNH}. By Proposition \ref{prop2}, $E$ Hyp-reduces to
$=_2$. By Proposition \ref{prop1}, $=_2$ does not Hyp-reduce to $E$.

To prove the second statement, we let $E$ consist of exactly $n+1$
equivalence classes, such that only $n$ of them contain Hyp reals.
For each $i< n-1$, we define the $i$th equivalence class by taking
all $x\in\sim X$, such that $x(0)=i$. We take the $n$th class to
contain all $x\in\sim X$ with $x(0)\geq n-1$. And the $(n+1)$st
class is $X$.

For the proof of the third statement, consider an equivalence
relation with $n_1$ classes such that only $n_0$ of them contain Hyp
reals.
\end{proof}

\begin{thm}
There are incomparable Hyp equivalence relations between $=_1$ and
$=_2$.
\end{thm}

\begin{proof}
To prove the theorem, we consider the following equivalence
relations. Let $A$ and $B$ be as in Theorem \ref{MainThm}. We take
the equivalence relation $E_A$ with two equivalence classes $A,\sim
A$ and $E_B$ with two equivalence classes $B, \sim B$. Then $E_A$
and $E_B$ are Hyp-reducible to $=_2$. By the properties of $A$ and
$B$, the relations $E_A$ and $E_B$ are Hyp-incomparable, as
otherwise (using the Remark following Theorem \ref{MainThm}) we
would have a Hyp function which maps $A$ to $B$ or vice versa.
\end{proof}

\begin{thm}
The partial order of Hyp subsets of $\omega$ under inclusion can be
order-preservingly embedded into the structure of degrees of Hyp
equivalence relations between $=_1$ and $=_2$.
\end{thm}

\begin{proof}
Let $X$ be a Hyp subset of $\omega$. Define the corresponding
equivalence relation $E_X$ in the following way. We let $xE_Xy$ iff
both $x,y\in\bigcup_{i\in X}A_i$ or both $x,y\in\sim\bigcup_{i\in
X}A_i$, where $A_0,A_1,\ldots$ are the sets constructed in Theorem
\ref{MainInf}. We check that $X\subseteq Y\Longleftrightarrow
E_X\leq_HE_Y$.

Suppose $X\subseteq Y$. For every $i\in X$ we send $A_i$ into
itself. We send $\sim~\bigcup_{i\in X}A_i$ into a single Hyp real
chosen in $\sim\bigcup_{i\in Y}A_i$. Therefore $E_X\leq_HE_Y$.

Now suppose $X\nsubseteq Y$ but $E_X\leq_H E_Y$ via a Hyp function
$F:\baire\to\baire$. Note that neither $\bigcup_{i\in X}A_i$ nor
$\bigcup_{i\in Y}A_i$ contain Hyp reals. Thus $F$ sends
$\sim\bigcup_{i\in X}A_i$ to $\sim\bigcup_{i\in Y}A_i$ and
$\bigcup_{i\in X}A_i$ to $\bigcup_{i\in Y}A_i$. Choose an $i_0\in
X\setminus Y$. Then $F[A_{i_0}]\subseteq\bigcup_{i\in
Y}A_i\subseteq\bigcup_{i\neq i_0}A_i$, contradicting the properties
of the sequence $A_0,A_1\ldots$
\end{proof}

\begin{cor}
 \begin{enumerate}
  \item There are infinite antichains between $=_1$ and $=_2$.
  \item There are infinite descending chains between $=_1$ and $=_2$.
  \item There are infinite ascending chains between $=_1$ and $=_2$.
 \end{enumerate}
\end{cor}

The same proof shows:
\begin{cor}
For any $1\leq n_0 < n_1\leq\omega$ there is an embedding of
$\mathcal{P}(\omega)\cap $Hyp into the structure of degrees of Hyp
equivalence relations that are above $=_{n_0}$, below $=_{n_1}$ and
incomparable with each $=_n$ for $n_0<n<n_1$.
\end{cor}

\subsection{Hyp Equivalence Relations between $=_{\omega}$ and $=_{\mathcal{P}(\omega)}$}

Let $=_{\mathcal{P}(\omega)}$ denote the Hyp-degree of the
equivalence relation of $=$ on $\mathcal{P}(\omega)$. By Proposition
\ref{prop2} and Silver's dichotomy \cite{gao09}, every Hyp
equivalence relation $E$ is either Hyp reducible to $=_{\omega}$, or
$=_{\mathcal{P}(\omega)}$ is Borel reducible to $E$. In
\S\ref{optsec} we will show that ``Borel reducible'' can be taken to
be ``Hyp in Kleene's $O$ reducible'', and that this is best
possible.



\begin{thm}
There exist Hyp-incomparable Hyp equivalence relations between
$=_{\omega}$ and $=_{\mathcal{P}(\omega)}$.
\end{thm}

\begin{proof}
Suppose that $A$ and $B$ are the $\Pi^0_1$ sets from Theorem
\ref{MainThm}: they contain no Hyp reals and there is no Hyp
function $F$ such that $F[A]\subseteq B$ or $F[B]\subseteq A$.

Now consider the equivalence relations $E_{A}$ and $E_{B}$:
$$
xE_{A}y \iff [(x\in A \wedge x=y) \vee
               (x,y\notin A \wedge x(0)=y(0))]
$$
and similarly for $E_B$ with $B$ replacing $A$.

By sending $n$ to the real $(n,0,0,\ldots)$ we get a Hyp reduction
$=_{\omega}$ to $E_{A}$ and $E_{B}$. Also $E_{A}$ (resp. $E_B$)
Hyp-reduces to $=_{\mathcal{P}(\omega)}$ via the map $G(x)=x$ if $x$
belongs to $A$ (resp. $B$), $G(x)=(x(0),0,0,\ldots)$ for $x\notin A$
(resp. $x\notin B$).

There is no Hyp reduction of $E_{A}$ to $E_{B}$. Indeed, suppose
that $F$ were such a reduction and let $C$ be the preimage under $F$
of $\sim B$. As $\sim B$ is $\Sigma^0_1$, $C$ is Hyp and therefore
$A\cap C$ is also Hyp. But $A\cap C$ must be countable as $F$ is a
reduction. So by Fact \ref{fct0}, part 2, if $A\cap C$ were nonempty
it would have a Hyp element, contradicting the fact that $A$ has no
Hyp elements. Therefore $F$ maps $A$ into $B$, which is impossible
by the choice of $A$ and $B$.
\end{proof}

\begin{thm}\label{thm11}
The partial order of Hyp subsets of $\omega$ under inclusion can be
embedded into the structure of degrees of Hyp equivalence relations
between $=_{\omega}$ and $=_{\mathcal{P}(\omega)}$.
\end{thm}

\begin{proof}
Let $A_0,A_1,\ldots$ be the $\Pi^0_1$ sets from Theorem
\ref{MainInf}. For every Hyp set $X\subseteq\omega$ consider the
equivalence relation of the form
$$
xE_{X}y \iff [(x\in \bigcup_{i\in X}A_i \text{ and } x=y) \text{ or }
               (x,y\notin \bigcup_{i\in X}A_i \text{ and }x(0)=y(0))].
$$
Then $=_{\omega}\leq_H E_X\leq_H =_{\mathcal{P}(\omega)}$. Suppose
$X\subseteq Y$. Then $E_{X}$ Hyp-reduces to $E_Y$ via the map
$G(x)=x$ if $x\in\bigcup_{i\in X}A_i$, $G(x)=(x(0),0,0,\ldots)$ for
$x\notin \bigcup_{i\in X}A_i$.

Suppose $X\nsubseteq Y$ but $E_X\leq_H E_Y$ via a Hyp function $F$.
Pick $i_0\in X\setminus Y$. As before, we consider the set
$$
A_{i_0}\cap F^{-1}(\sim\bigcup_{j\in Y}A_j).
$$
Then this is a countable Hyp set. If it is non-empty then it
contains a Hyp real, contradicting the definition of $A_{i_0}$.
Therefore we get $F[A_{i_0}]\subseteq\bigcup_{j\in Y}A_j\subseteq
\bigcup_{j\neq i_0}A_j$, contradiction.
\end{proof}

\begin{cor}
There are infinite chains and antichains between $=_{\omega}$ and
$=_{\mathcal{P}(\omega)}$.
\end{cor}

\begin{cor}
For any finite $n_0\geq 1$, the partial order of Hyp subsets of $\omega$
under inclusion can be embedded into the structure of degrees of Hyp
equivalence relations between $=_{n_0}$ and
$=_{\mathcal{P}(\omega)}$ but incomparable with $=_n$ for
$n_0<n\leq\omega$.
\end{cor}

\begin{proof}
For every Hyp $X\subseteq\omega$, consider the equivalence relation
of the form
\begin{align*}
 xE^{n_0}_{X}y \iff & x\in \bigcup_{i\in X}A_i\wedge x=y\vee\\
                    & x,y\notin \bigcup_{i\in X}A_i\wedge(x(0)=y(0)<n_0-1\vee x(0),y(0)\geq n_0-1).
\end{align*}
Then $E^{n_0}_{X}$ has exactly $n_0$ equivalence classes with Hyp
reals. Therefore $=_{n_0}\leq_H E^{n_0}_X$ and for
$n_0<n\leq\omega$, the equivalence relation $=_n$ is incomparable
with $E^{n_0}_{X}$.
\end{proof}

\subsection{Hyp Equivalence Relations between $=_{\mathcal{P}(\omega)}$ and $E_0$}

It was shown in \cite{hakelou90} that any Hyp equivalence relation
is either Hyp reducible to $=_{\mathcal{P}(\omega)}$, or $E_0$ is
Borel reducible to it. In \S\ref{optsec} we will show that ``Borel''
can be taken to be ``Hyp in Kleene's ${O}$'', and that this is best
possible.

\begin{thm}\label{thm12}
There exist Hyp-incomparable Hyp equivalence relations between
$=_{\mathcal{P}(\omega)}$ and $E_0$.
\end{thm}

\begin{proof}
Let $A$ and $B$ be the Hyp sets from Theorem \ref{MainThm}, such
that for no Hyp function $F$ do we have $F[A]\subseteq B$ nor
$F[B]\subseteq A$.

Define two Hyp equivalence relations $E_{A}$ and $E_{B}$ on $\baire\times
2^\omega$ by
$$
(x,y) E_{A} (x',y')\iff x=x'\wedge[(x\notin A)\vee(x\in A\wedge yE_0 y')],
$$
and similarly for $E_B$ with $B$ replacing $A$.

Suppose $F:\baire\times 2^\omega\to\baire\times 2^\omega$ is a
Hyp-reduction of $E_{A}$ to $E_{B}$. Define $F'(x,y)=z\iff (\exists
w) F(x,y)=(z,w)$. Note that $F'$ is constant on $E_{A}$ classes.
Define a function $h:\baire\to \baire$ by
\begin{align*}
h(x)=z&\iff \{y\in 2^\omega: g'(x,y)=z\}\text{ is non-meagre}\\
      (&\iff \{y\in 2^\omega: g'(x,y)=z\}\text{ is comeagre.})
\end{align*}
By Facts \ref{fct1} and \ref{fct2}, $h$ is an everywhere defined Hyp
function. Suppose $x\in A$. Then for a comeagre set $C\subseteq
2^\omega$ we have $F'(x,y)=h(x)$ for all $y\in C$. We claim that
$h(x)\in B$. Indeed, otherwise the set $\{x\}\times C$ is mapped by
$g$ into a single $E_{B}$ class, contradicting that all
$E_{A}|\{x\}\times 2^\omega$ classes are meagre in $\{x\}\times
2^\omega$ (in fact, they are countable).

Thus $h$ is a Hyp function with $h[A]\subseteq B$, contradicting the
properties of $A$ and $B$.
\end{proof}

\begin{thm}
The partial order of Hyp subsets of $\omega$ can be embedded into
the structure of Hyp equivalence relations between
$=_{\mathcal{P}(\omega)}$ and $E_0$.
\end{thm}
\begin{proof}
Let $A_0, A_1,\ldots$ be the sequence from Theorem \ref{MainInf}.
For every Hyp set $X$ we define a Hyp equivalence relation $E_X$ on
$\baire\times 2^\omega$ in the following way:
$$
(x,y) E_{X} (x',y')\iff x=x'\wedge[(x\notin\bigcup_{i\in X} A_i)\vee(x\in\bigcup_{i\in X} A_i\wedge yE_0 y')].
$$
Then the theorem follows from an argument similar to that in the
proof of Theorems \ref{thm11} and \ref{thm12}.
\end{proof}

\begin{thm}
For any $n_0\leq\omega$ the partial order of Hyp subsets of $\omega$
can be embedded into the structure of degrees of Hyp equivalence
relations between $=_{n_0}$ and $E_0$, but incomparable with $=_n$
for $n_0<n\leq\omega$ and incomparable with
$=_{\mathcal{P}(\omega)}$.
\end{thm}

\begin{proof}
Let $A_0,A_1,\ldots$ be the sequence of Hyp sets from Theorem
\ref{MainInf}. For a Hyp set $X\subseteq\omega$, define an
equivalence relation $E_{X}^{n_0}$ on $\baire\times 2^\omega$ by
\begin{align*}
(x,y) E^{n_0}_{X} (x',y')\iff & [x,x'\notin\bigcup_{i\in X} A_i\wedge \\
&(x(0)=x'(0)<n_0-1\vee x_0,x'_0\geq n_0-1)]\vee\\
& [x,x'\in\bigcup_{i\in X} A_i\wedge x=x'\wedge yE_0 y'].
\end{align*}
The relation $E^{n_0}_X$ is Hyp. Clearly it is below $E_0$. It has
only $n_0$ equivalence classes with Hyp reals, thus it is above
$=_{n_0}$ and incomparable with $=_n$ for $n_0<n<\omega$ and with
$=_{\mathcal{P}(\omega)}$.
\end{proof}

\section{Category notions in the Gandy-Harrington topology}

The spaces under consideration in this section will be of the form
$(\omega^\omega)^n$, $1\leq n\leq\omega$. Baire space,
$\baire=\omega^\omega$, is a Polish space in the product topology,
and thus so is $\baire^n$ for all $n\leq\omega$. We will call this
the ``usual'' topology on $\baire^n$. We consider two other topologies on $\baire^n$:
\begin{enumerate}[(1)]
\item The {\it Gandy-Harrington} topology, which is generated by
the (lightface) $\Sigma^1_1$ subsets of $\baire^n$. This topology
will be denoted $\tau_n$ if $n>1$, or simply by $\tau$ if $n=1$.
\item The product topology $\tau^n$ on $\baire^n$, when we equip
$\baire$ with the Gandy-Harrington topology.
\end{enumerate}

These topologies are all different: The usual topology is weaker
than $\tau^n$, which again is weaker than $\tau_n$, if $n>1$.

The purpose of this section is to examine the effectiveness of
category notions in the Gandy-Harrington topology. For instance,
if we consider a $\Sigma^1_1$ set $A\subseteq \baire^2$, we would
like to know the complexity of the set
$$
\{x\in\baire: A_x\text{ is not meagre in }\tau\}.
$$
We would also like to know how effective we can reasonably expect
a winning strategy in the Banach-Mazur game to be, or how
effective player II's winning strategy in the Choquet game in
$(\baire,\tau)$ is. Our analysis is entirely parallel to that
found in \cite{kechris73}, where the same questions were analyzed
for the usual topology on $\baire$.

It is important to note that the category of a set may change when
changing between these topologies. For instance, a $\Sigma^1_1$
singleton $\{x\}\subseteq \baire$ is open in the Gandy-Harrington
topology, but closed and meagre in the usual topology. The set
$$
A=\{(x,x): \{x\} \text{ is not a $\Sigma^1_1$ singleton}\}
$$
is meagre in $\tau^2$ since every section $A_x$ is meagre, but it
is open in $\tau_2$ since it is $\Sigma^1_1$. However, all
$\bSigma^1_1$ subsets of $\baire^n$ have the property of Baire in
the topologies $\tau^n$ and $\tau_n$. This follows from
\cite[Theorem 21.8]{kechris95}.

\subsection{Basic computations}

Fix $n$ and let $\sigma$ be either $\tau_n$ or $\tau^n$. Then we
form the finite levels of the Borel hierarchy:
$\bSigma^0_1[\sigma]$ consists of the $\sigma$-open subsets of
$\baire^n$, and in general $\bSigma^0_{k+1}[\sigma]$ consists of
countable unions of sets from $\bPi^0_k[\sigma]$, which itself
consists of the complements of sets in $\bSigma^0_{k}[\sigma]$.

Let $A\subseteq\omega\times\baire^n$ be universal for
$\Sigma^1_1$. Then
$$
G_1=\{(x,y)\in\baire\times\baire^n: (\exists k) x(k)>0\wedge y\in
A_k\}
$$
is a $\Sigma^1_1$ set which is universal for
$\bSigma^0_1[\tau_n]$. The set $\baire\times\baire^n\setminus G_1$
is $\Pi^1_1$ and universal for $\bPi^0_1[\tau_n]$, and
$$
G_2=\{(x,y)\in \baire^\omega\times\baire^n: (\exists k) y\in
(\baire\times\baire^n\setminus G_1)_{x(k)}\}
$$
is $\Pi^1_1$ and universal for $\bSigma^0_2[\tau_n]$, and we can
continue in this way to find universal sets for
$\bSigma^0_k[\tau_n]$ that are $\Sigma^1_1$ when $k$ is odd and
$\Pi^1_1$ when $k$ is even. A similar analysis applies to
$\tau^n$.

\begin{prop}
Let $\sigma$ be either $\tau_n$ or $\tau^n$, and let
$A\subseteq\baire\times\baire^n$ be a $\Sigma^1_1$ or $\Pi^1_1$
universal set for $\bSigma^0_k[\sigma]$, depending on if $k$ is
odd or even. Then
$$
\{x\in\baire: A_x\text{ is not } \sigma\text{-meagre}\}
$$
is $\Delta^1_1(O)$, where $O$ denotes Kleene's $O$.\label{finProp}
\end{prop}

\begin{proof}
Let $\sigma=\tau_n$. If $A\in\bSigma^0_1[\sigma]$ then there is a
sequence $T_l$ of recursive trees on $\omega^{n+1}$ such that
$$
A=\bigcup_{l\in\omega} p[T_l],
$$
where $p[T_l]$ is the projection of the set $[T_l]$ of infinite branches through $T_l$. Now
$$
A\text{ is not meagre }\iff (\exists l) [T_l]\neq\emptyset.
$$
is clearly arithmetic in the sequence $(T_l)$ and Kleene's $O$.

If $A\in\bSigma^0_{k+1}[\sigma]$, find a sequence
$B_l\in\bPi^0_k[\sigma]$ such that
$$
A=\bigcup_{l\in\omega} B_l.
$$
Then
\begin{align*}
A\text{ is not meagre }&\iff (\exists l) B_l \text{ is not meagre}\\
&\iff (\exists l)(\exists T\text{ recursive}) p[T]\setminus B_l \text{ is meagre}\\
&\iff(\exists l)(\exists T\text{ recursive}) \neg (p[T]\setminus B_l \text{ is not meagre})\\
\end{align*}
which is arithmetic in $O$ and the sequence $(B_l)$ by the
inductive hypothesis. The proof of the case $\sigma=\tau^n$ is
similar.
\end{proof}

Our next goal is to prove the following:

\begin{prop}
Let $\sigma$ be $\tau_n$ or $\tau^n$, and let
$A\subseteq\baire\times\baire^n$ be a $\Sigma^1_1$ set universal
for $\bSigma^1_1(\baire^n)$. Then
$$
\{x\in\baire: A_x \text{ is not } \sigma\text{-meagre}\}
$$
is $\Sigma^1_1(O)$.\label{sigma11prop}
\end{prop}

Before proving this we need a generalization of Proposition 1.5.2
in \cite{kechris73}. Let $(X,\sigma)$ be a 2nd countable
topological space and let $\mathcal U$ be a countable basis for
the topology.

A function $f:\mathcal U\to\omega^{<\omega}$ is called {\it
$\mathcal U$-monotone} if
$$
(\forall U,V\in\mathcal U) U\subseteq V\implies f(V)\subseteq
f(U).
$$
For $x\in\baire$ we define
$$
(\lim_{\mathcal U} f)(x)=y\iff (\forall k)(\exists U\in\mathcal U)
x\in U\wedge \lh(f(U))\geq k\wedge f(U)\subseteq y.
$$
The set
$$
\{x\in\baire: (\exists y)(\lim_{\mathcal U} f) (x)=y\}
$$
is $G_\delta$ in the topology $\sigma$, and $\lim_{\mathcal U} f$
defines a function on this set. With these definitions we have the
following analogue of \cite[Proposition 1.5.2]{kechris73}:

\newtheorem{lemma}{Lemma}
\begin{lemma}[Folklore]
Let $(X,\sigma)$ be a 2nd countable topological space and let
$\mathcal U$ be a countable basis for the topology $\sigma$. Then:

\begin{enumerate}[\rm (1)]

\item If $Y\subseteq X$ is a $G_\delta$ set and $\bar
f:Y\to\omega^\omega$ is continuous w.r.t. the usual topology on
$\omega^\omega$ and $\sigma$ on $Y$, then there is a $\mathcal
U$-monotone function $f:\mathcal U\to\omega^{<\omega}$ such that
$\bar f=\lim_{\mathcal U} f$.

\item If $f:\mathcal U\to\omega^{<\omega}$ is $\mathcal
U$-monotone then $\lim_{\mathcal U} f$ is continuous on its domain
(taking $\omega^\omega$ with the usual topology and
$\dom(\lim_{\mathcal U} f)$ with the topology induced by
$\sigma$), and $\dom(\lim_{\mathcal U} f)$ is a $G_\delta$ set in
the topology $\sigma$.

\end{enumerate}
\label{folklorelemma}
\end{lemma}

\begin{proof}
(2) is clear from the definition. For (1), let
$Y=\bigcap_{n\in\omega} W_n$, where the $W_n$ are open sets. Let
$(U_n)$ enumerate $\mathcal U$. We can assume that
$W_{n+1}\subseteq W_n$. Define $f:\mathcal U\to\omega^{<\omega}$
by letting $f(U_k)$ be the longest sequence $s$ such that
$$
(\forall l)(U_l\subseteq U_k\implies \bar f(U_l)\subseteq
N_s)\wedge \lh(s)\leq\min\{k,\max\{n:U_k\subseteq W_n\}\}.
$$
(Here $N_s$ denotes the basic open neighborhood determined by $s$,
i.e.
$$
N_s=\{x\in\baire: s\subseteq x\}.)
$$
Since $\bar f$ is continuous it follows that if $x\in Y$ then for
all $k,n\geq 0$ we can find $U_m\subseteq W_n$ such that $x\in
U_m$ and $f(U_m)\subseteq N_s$ for some sequence $s$ of length at
least $k$. Thus $x\in\dom(\lim_{\mathcal U}f)$ and clearly $\bar
f(x)=(\lim_{\mathcal U} f)(x)$. On the other hand, if $x\notin Y$
then there is $n$ such that $x\notin W_n$. Thus $\lh(f(U_k))\leq
n$ for all $k\in\omega$, and so $x\notin\dom(\lim_{\mathcal U}
f)$.
\end{proof}

We now turn to the proof of Proposition \ref{sigma11prop}. Recall that the set
$$
X_{\low}=\{x\in\omega^\omega: \omega_1^{x}=\omega_1^{\ck}\}
$$
is $\Sigma^1_1$ and furthermore that it is dense in the
Gandy-Harrington topology, see e.g. Appendix A of \cite{gao09}.

\begin{proof}[Proof of Proposition \ref{sigma11prop}]
The proof follows the general lines of \cite[Theorem
2.2.5]{kechris73}. For simplicity we consider the case $n=1$, i.e.
$\sigma=\tau$. Let $A\subseteq \baire$ be $\bSigma^1_1$ and not
meagre, and let $T$ be a tree on $\omega\times\omega$ such that
$p[T]=A$. Then by the Jankov-von Neumann uniformization theorem
\cite[18.1]{kechris95} we may find a $\mathcal B(\bSigma^1_1)$
uniformising function $\bar f:A\to \baire$ such that
$$
(\forall x\in A) (x,\bar f(x))\in [T].
$$
(Here $\mathcal B(\bSigma^1_1)$ denotes the $\sigma$-algebra
generated by the $\bSigma^1_1$ sets.) Since every $\bSigma^1_1$
set has the Baire Property in $\tau$ it follows that the function
$\bar f$ is Baire measurable when $\dom(\bar f)$ is given the
topology $\tau$ and $\codom(\bar f)$ is given the usual topology.
Since $A$ has the Baire Property in $\tau$ we may find a
$\tau$-$G_\delta$ set $A'\subseteq A$ such that
\begin{enumerate}[(a)]
\item $A\setminus A'$ is $\tau$-meagre, \item $\bar f|A'$ is
continuous (w.r.t. $\tau$ on $\dom(\bar f)$ and the usual topology
in $\codom(\bar f)$.) \item $A'\subseteq X_{\low}$
\end{enumerate}
Now let $B\subseteq\omega\times\baire$ be $\Sigma^1_1$ such that
$$
C\subseteq X_{\low} \text{ is } \Sigma^1_1\iff (\exists n) C=B_n.
$$
Let $\mathcal U=\{B_n: B_n\neq\emptyset\}$. Then by Lemma
\ref{folklorelemma} we can find a monotone $f:\mathcal
U\to\omega^{<\omega}$ such that $\bar f|A'=\lim_{\mathcal U} f$
and
$$
(\forall s)(\forall n)(s|\lh(f(B_n\cap N_s)),f(B_n\cap
N_s)|\lh(s))\in T.
$$
Thus
\begin{align*}
&A\text{ is not meagre}\iff\\
&(\exists f:\omega\to\omega^{<\omega}) ((\forall n) (B_n=\emptyset\implies f(n)=\emptyset)\wedge\\
&(\forall m)(\forall n)(B_m\subseteq B_n\wedge
B_m\neq\emptyset\implies f(n)\subseteq
f(m))\wedge\\
& (\forall s\in\omega^{<\omega})(\forall
n)(s|\lh(f(n)),f(n)|\lh(s))\in
T\wedge \\
&\dom(\lim_{\mathcal U} f) \text{ is not meagre}).
\end{align*}
where above, $\lim_{\mathcal U} f$ has the natural meaning if we
think of $f$ as being defined on $\mathcal U$, not on the indices
of elements of $\mathcal U$.

If $f:\omega\to\omega^{<\omega}$ is (a code for a) monotone
function then
$$
x\in\dom(\lim_{\mathcal U} f)\iff(\forall k)(\exists n) x\in
B_n\wedge \lh(f(n))\geq k,
$$
so ``$\dom(\lim_{\mathcal U} f) \text{ is not meagre}$'' is
$\Delta^1_1(O,f)$ uniformly in $f$ by Proposition \ref{finProp}.

The proof is finished by noting that the statement ``$B_m\subseteq
B_n$'' may be replaced by the statement
$$
\neg(B_m\setminus B_n  \text{ is not meagre}).
$$
To see this, note that by \cite[Theorem A.1.6]{gao09} we have that
if $D\subseteq X_{\low}$ is $\Sigma^1_1$ then $D$ is $\tau$-clopen
in $X_{\low}$. Thus $B_m\setminus B_n=\emptyset$ iff $B_m\setminus
B_n$ is meagre. Since by Proposition \ref{finProp} the statement
``$B_m\setminus B_n$ is not meagre'' is $\Delta^1_1(O)$, this
finishes the proof.
\end{proof}

\subsection{The Choquet and Banach-Mazur games}
Let $\sigma=\tau_n$ or $\sigma=\tau^n$. Recall the {\it strong
Choquet game} $G_{(\baire^n,\sigma)}$:
$$
\begin{array}{rlllll}
\mathrm{I} & x_0, U_0 & \   & x_1, U_1\\
  &          &     &          &      & \cdots\\
\mathrm{II} & \       & V_0 &\         &  V_1

\end{array}
$$
Players I and II take turns playing. The $i$th move for Player I
consists of a basic open set $U_i$ and a point $x_i\in U_i$.
Player II must respond by playing a basic open set $V_i\subseteq
U_i$ such that $x_i\in V_i$. Then Player I is required to
respond with $x_{i+1}$ and $U_{i+1}$ such that $x_{i+1}\in
U_{i+1}\subseteq V_i$. Player II wins iff
$$
\bigcap_{n\in\omega} U_n=\bigcap_{n\in\omega} V_n\neq\emptyset.
$$
It is well-known that II has a winning strategy in the strong
Choquet game in $(\baire,\tau)$, see e.g. \cite[Theorem
4.1.5]{gao09}. Moreover, the winning strategy for II described in
the proof there is $\Delta^1_1$ in the codes. From this we easily
get:
\begin{cor}
Let $n\geq 1$ and let $\sigma=\tau_n$ or $\sigma=\tau^n$. Then
{\rm II} has a winning strategy in $G_{(\baire^n,\sigma)}$ which
is $\Delta^1_1$ in the codes.
\end{cor}

What about the Banach-Mazur game in $(\baire^n,\sigma)$? Recall
that the Banach-Mazur game $G^{**}_{\sigma}(A)$, where
$A\subseteq\baire^n$ is non-empty, is played as follows: Players I
and II take turns playing basic open sets $U_i$ and $V_i$,
$$
\begin{array}{rlllll}
\mathrm{I} & U_0 & \   & U_1\\
  &          &     &          &      & \cdots\\
\mathrm{II} & \       & V_0 &\         &  V_1
\end{array}
$$
and the players are required to maintain that $U_i\supseteq
V_i\supseteq U_{i+1}$ for all $i\geq 0$. II wins iff
$$
\bigcap_{n\in\omega} U_n=\bigcap_{n\in\omega} V_n\subseteq A.
$$
It is well-known (see e.g. \cite[8.33]{kechris95}) that
$A\subseteq \baire^n$ is comeagre (in $\sigma$) if and only if II
has a winning strategy in $G_\sigma^{**}(A)$. By (ii) of
\cite[8.33]{kechris95}, it also follows that $A$ is meagre in a
non-empty open set if and only if I has a winning strategy.

In the case $\sigma=\tau_n$ any $\Delta^1_1$ set
$A\subseteq\baire^n$ is of course $\sigma$-clopen, and so if
$\baire^n\setminus A\neq\emptyset$ then I clearly wins simply by
playing $\baire^n\setminus A$ in the first move. For
$\sigma=\tau^n$ the situation is more complicated:

\begin{prop}\label{propA3}
Let $n\geq 2$. If $A\subseteq \baire^n$ is $\Delta^1_1$ then there
is a winning strategy in $G^{**}_{\tau^n}(A)$ which is
$\Delta^1_1(O)$ in the codes.
\end{prop}

\begin{proof}
The proof is a variation of \cite[Theorem 4.2.1]{kechris73}. For
notational simplicity, we deal with the case $n=2$. Moreover,
following \cite[Definition 8.25]{kechris95} we will use the
notation
$$
U\forces A,
$$
where $U$ is a basic open set and $A$ some subset, to mean that
$A$ is comeagre in $U$, i.e. $U\setminus A$ is meagre. Finally, we
fix a $\Delta^1_1$ winning strategy for II in the strong Choquet
game in $(\baire^2,\tau^2)$. Since there plainly is a danger of
confusion here, we will refer to the players of the {\it strong
Choquet game} as $\I_C$ and $\II_C$. I and II then refers to the
players in the {\it Banach-Mazur game}.

Without loss of generality, assume that I wins
$G_{\tau^2}^{**}(A)$, i.e. $\baire^2\setminus A$ is comeagre in
some open set. We will describe a $\Delta^1_1(O)$ winning strategy
for I. Player I will be aided by playing (as $\I_C$) a strong
Choquet game concurrently with the Banach-Mazur game.
Schematically:

$$
\begin{array}{lrlllllllll}
\mathrm{I}_C & :  & x_0, B_{k_0} & \       & \        & \        & x_2, B_{k_2} &      &  \\
               &         &         &          &          &         &      &           &  & &\cdots\\
\mathrm{II}_C&: & \       & B_{n_0} & \        &          &         & B_{n_2}          &  \\
\\

\mathrm{I}&:    & \       & \       & B_{n_0}  &  \       &         &      &  B_{n_2}  &   & \\
               &         &         &          &          &         &      &           &   & &\cdots\\
\mathrm{II}&:   &         &         &          & B_{n_1}  &         &      &           &  B_{n_3} &  \\

\end{array}
$$

Fix a $\Delta^1_1$-scale $\{\varphi_m\}_{m\in\omega}$ on
$\baire^2\setminus A$. For $x\in\baire^2$ let
$$
\psi_m(x)=\langle\varphi_0(x),x(0),\cdots,\varphi_m(x),x(m)\rangle,
$$
where as in \cite{kechris73},
$\langle\gamma_0,\ldots\gamma_m\rangle$ is the rank of
$(\gamma_0\ldots,\gamma_m)$ in the lexicographic order on
$\mathbb{ON}^{<\omega}$. Note that $X_{\low}^2$ is open and dense in
$(\baire,\tau^2)$. Let $B\subseteq\omega\times\baire^2$ be a
$\Sigma^1_1$ parametrization of
$$
\{C_0\times C_1:C_0,C_1\in\Sigma^1_1(\baire), C_0,C_1\subseteq
X_{\low}\}.
$$
First find $B_{k_0}$ where $k_0$ is least such that
$B_{k_0}\neq\emptyset$ and $B_{k_0}\forces \baire^2\setminus A$,
and let $x_0\in B_{k_0}$ be computable in $O$. Let $B_{n_0}$ be
the response of $\II_C$ according to the fixed winning strategy in
the strong Choquet game when $\I_C$ plays $x_0, B_{k_0}$. I's
first move in the Banach-Mazur game is then $B_{n_0}$. Suppose II
responds by playing $B_{n_1}$. Let $s_1\in \omega^{<\omega}$ be
the least sequence of length 1 such that $B_{n_1}\cap N_{s_1}$ is
not meagre. Now
\begin{align*}
A_2=&\{x: x\in \baire^2\setminus A\wedge x\in B_{n_1}\cap N_{s_1}\wedge\\
&\{y: y\in \baire^2\setminus A\wedge y\in B_{n_1}\cap
N_{s_1}\wedge
\psi_1(x)=\psi_1(y)\}\text{ is not meagre }\wedge\\
&\{y:y\in \baire^2\setminus A\wedge x\in B_{n_1}\cap
N_{s_1}\wedge\psi_1(y)<\psi_1(x)\}\text{ is meagre}\}\}
\end{align*}
is non-meagre, $\Delta^1_1(O)$ and $A_2\subseteq N_{s_1}\cap
\baire^2\setminus A$.

Now let $k_2$ be least such that $B_{k_2}\neq\emptyset$,
$B_{k_2}\subseteq B_{n_1}$ and $B_{k_2}\forces A_2$. We may find
$k_2$ in a $\Delta^1_1(O)$ way since, as in the proof of
Proposition \ref{sigma11prop}, $B_{k_2}\subseteq B_{n_1}$ may be expressed by
saying that $B_{k_2}\setminus B_{n_1}$ is meagre since we work on
$X_{\low}^2$. Let $x_2\in B_{k_2}$ be computable in $O$ and let
$\I_C$ play  $x_2, B_{k_2}$ in the strong Choquet game. $\II_C$
responds with $B_{n_2}$. I plays $B_{n_2}$ in the Banach-Mazur
game.

Suppose II responds by playing $B_{n_3}$. Let $s_3$ be the least
sequence of length 3 such that $s_1\subseteq s_3$ and $N_{s_3}\cap
B_{n_3}$ is not meagre. We let
\begin{align*}
A_4=&\{x: x\in A_2\wedge x\in B_{n_3}\cap N_{s_3}\wedge \{y: y\in A_2\wedge y\in B_{n_3}\cap N_{s_3}\wedge\\
&\psi_3(x)=\psi_3(y)\}\text{ is not meagre }\wedge\\
&\{y:y\in A_2\wedge y\in B_{n_3}\cap
N_{s_3}\wedge\psi_3(y)<\psi_3(x)\}\text{ is meagre}\}\}.
\end{align*}
Then $A_4\subseteq N_{s_3}\cap A_2$ and $A_4$ is $\Delta^1_1(O)$
and non-meagre, and so we let $\I_C$ play $x_4, B_{k_4}$ in the
strong Choquet game, where $k_4$ is least such that $\emptyset\neq
B_{k_4}\subseteq B_{n_3}$, $B_{k_4}\forces A_4$ and $x_4\in
B_{k_4}$ is computable in $O$. $\II_C$ responds with $B_{n_4}$,
and I plays $B_{n_4}$ in the Banach Mazur game, and so on. At the
end of this run of the Banach-Mazur game we have produced a
sequence of sets $B_{n_i}$ and a real $\alpha=\bigcup_{i\in\omega}
s_{2i+1}$.

Note that since $\II_C$ wins the strong Choquet game we have that
$$
\bigcap_{i\in\omega} B_{n_i}\neq\emptyset.
$$
Clearly, it must then be the case that
$$
\bigcap_{i\in\omega} B_{n_i}=\{\alpha\}.
$$
We claim that $\alpha\notin A$. For this, note that by
construction we can find a sequence $(x_i)$ such that $x_i\in
A_{2i}$ and $s_{2i-1}\subseteq x_i$. For this sequence it holds
for all $m$ that $\psi_{2m+1}(x_i)$ is constant for $i\geq m$.
Thus $\varphi_m(x_i)$ is eventually constant, and since $x_i\to
\alpha$ it follows by the properties of a scale that $\alpha\in
\baire^2\setminus A$. Thus I wins this run of the game.
\end{proof}

\noindent \textbf{Remark.} The previous proof relativizes to a parameter in the
following way: If $A\subseteq \baire^n$ is $\Delta^1_1(z)$ for
some real $z$ then one of the players has a $\Delta^1_1(O,z)$
winning strategy in the game $G^{**}_{\tau^n}(A)$. It is easy to
see that the same proof also goes through for
$G^{**}_{\tau_n}(A)$. Thus we have

\begin{cor}
Let $n\geq 1$ and let $\sigma=\tau^n$ or $\sigma=\tau_n$. If
$A\subseteq\baire^n$ is $\Delta^1_1(z)$ then there is a winning
strategy in $G^{**}_{\sigma}(A)$ which is $\Delta^1_1(O,z)$ in the
codes. In particular, if $A$ is $\Delta^1_1(O)$ then there is a
$\Delta^1_1(O)$ winning strategy in $G^{**}_\sigma(A)$.
\end{cor}

From this we get:

\begin{cor}
Let $n\geq 1$ and $\sigma=\tau^n$ or $\sigma=\tau_n$. Suppose
$A\subseteq\baire^n$ is $\Delta^1_1(O)$ and comeagre. Then there
is $C\subseteq\omega\times\baire^n$ such that
\begin{enumerate}[$(1)$]
\item $C$ is $\Sigma^1_1(O)$, \item For all $n\in\omega$, $C_n$ is
$\sigma$-open and dense \item $\bigcap_{n\in\omega} C_n\subseteq
A$.
\end{enumerate}
\label{effcomeagre}
\end{cor}

\begin{proof}
As in the proof of Corollary 4.2.4 in \cite{kechris73}, we need
only note that the construction given in the proof of Theorem 6.1
in \cite{oxtoby96} produce the desired set $C$. To see this, fix a
$\Sigma^1_1$ set $B\subseteq\omega\times\baire^n$ such that
$$
U\subseteq\baire^n \text{ is a basic $\sigma$-open set}\iff
(\exists n) U=B_n.
$$
If we use the strategy described in Proposition \ref{propA3} above in the
proof of Theorem 6.1 in \cite{oxtoby96}, we will obtain a sequence
$W_n\subseteq\omega$ of $\Delta^1_1(O)$ sets (uniformly in $n$)
such that
$$
C_n=\bigcup_{k\in W_n} B_k
$$
is open dense and
$$
\bigcap C_n\subseteq A.
$$
Thus
$$
(n,x)\in C\iff (\exists k) k\in W_n\wedge x\in B_k
$$
gives a $\Sigma^1_1(O)$ definition of a set $B$ that is as
required.
\end{proof}

\section{Parameters in the basic dichotomy theorems}\label{optsec}

The results of the previous section show that complexity
computations involving category notions in the Gandy-Harrington
topology can be carried out using Kleene's $O$ as a parameter. In
this section we will use this (specifically, Corollary \ref{effcomeagre} above) to
show that the proofs of the Silver and Harrington-Kechris-Louveau
dichotomy Theorems produce reductions that are no worse than
$\Delta^1_1(O)$. We also show that this is in some sense the best
possible result we can hope for.

We start with Silver's dichotomy:

\begin{thm}[Silver's dichotomy]
Let $E$ be a Hyp equivalence relation on $\baire$. Then
either
$$
E\leq_{H} =_{\omega}
$$
or
$$
=_{\mathcal{P}(\omega)}\leq_{\Delta^1_1(O)} E.
$$
\label{silver}
\end{thm}
\begin{proof}
By Harrington's well-known proof of Silver's dichotomy (see e.g.
\cite[Theorem 32.1]{jech03} or \cite[Theorem 5.3.5]{gao09}),
either ({\it i}) every $E$-equivalence class contains a non-empty
Hyp set, or else ({\it ii}) there is a $\Sigma^1_1$ set
$H\subseteq\baire$ such that $E\cap H\times H$ is $\tau^2$-meagre
in $H\times H$. In case ({\it i}) by Proposition \ref{prop2} there is a
Hyp reduction of $E$ to $=_{\omega}$. We show that in case ({\it ii}) there is a $\Delta^1_1(O)$ reduction of
$=_{\mathcal{P}(\omega)}$ to $E$.

By Corollary \ref{effcomeagre} we may find a $\Sigma^1_1(O)$ set
$C\subseteq \omega\times\baire^2$ such that $C_n$ is $\tau^2$-open
dense in $H\times H$ and
$$
\bigcap_{n\in\omega} C_n\subseteq H\times H\setminus E.
$$
We may assume that $C_{n+1}\subseteq C_n$ for all $n$.
Harrington's proof (as presented in \cite{gao09} or \cite{jech03})
now produces a reduction of $=_{\mathcal{P}(\omega)}$ to $E$ which is
Hyp relative to the sequence $C_n$. To see this, fix a
Hyp winning strategy for II in the strong Choquet game
$G_{(\baire^2,\tau^2)}$. Then we may easily define a scheme
consisting of basic open sets $(U_s)_{s\in 2^{<\omega}}$,
$(V_s)_{s\in 2^{<\omega}}$ and points $(x_s)_{s\in 2^{<\omega}}$
such that $s\mapsto (U_s,V_s,x_s)$ is $\Delta_1^1(O)$ (in the codes) and the
following conditions hold:
\begin{enumerate}[(1)]
\item $U_\emptyset=V_\emptyset=H$ \item For each $s\in
2^{<\omega}$ the following is a play according to II's winning
strategy:
$$
\begin{array}{rlllllllll}
\mathrm{I}  & x_{s\restrict 1}, U_{s\restrict 1} &                  & x_{s\restrict 2}, U_{s\restrict 2}  &                  & \cdots  &        & x_s,U_s\\
   &                                   &                  &                                     &                  &         &        &         &  \\
\mathrm{II} &                                   & V_{s\restrict 1}
&                                     & V_{s\restrict 2} &
& \cdots &         & V_s
\end{array}
$$
\item $\diam(U_s)<2^{-\lh(s)}$ (with respect to the usual metric
on $\baire$). \item $U_{s\hat{\ }0}\times U_{s\hat{\ }1}\subseteq
C_{\lh(s)}$.
\end{enumerate}
If we define for $x\in 2^\omega$
$$
f(x)=y\iff y\in\bigcap_{n\in\omega} V_n
$$
then $f$ is a $\Delta^1_1(O)$ function and is easily seen to be a
reduction of $=_{\mathcal{P}(\omega)}$ to $E$.
\end{proof}

For the Glimm-Effros dichotomy
due to Harrington, Kechris and Louveau we have:

\begin{thm}[Harrington-Kechris-Louveau  \cite{hakelou90}]
Let $E$ be a Hyp equivalence relation on $\baire$. Then
either
$$
E\leq_{H} =_{\mathcal{P}(\omega)}
$$
or
$$
E_0\leq_{\Delta^1_1(O)} E.
$$\label{hkl}
\end{thm}
\begin{proof}
There are again two cases: (1) $E=E^*$, where $E^*$ is the closure
of $E$ in the topology $\tau^2$, or (2) $E\neq E^*$.

In the first case, it was observed in \cite{hakelou90}, p. 922,
that there is a Hyp reduction of $E$ to $=_{\mathcal{P}(\omega)}$.
So we only have to handle the 2nd case.

We will follow the exposition of the proof of the
Harrington-Kechris-Louveau Theorem given in \cite[\S 6.3]{gao09}.
Since $E\neq E^*$ the set
$$
X=\{x\in \baire: (\exists y) yE^* x\wedge \neg yEx\}.
$$
is non-empty and $\Sigma^1_1$. By \cite[Lemma 6.3.8]{gao09} $E$ is
dense and meagre in $X^2\cap E^*$. By Corollary \ref{effcomeagre},
we may find $C\subseteq\omega\times\baire^2$ such that $C_m$ is
$\tau^2$-open dense in $X^2$, $m\leq n\implies C_m\supseteq C_n$,
and
$$
\bigcap_{n\in\omega} C_n\subseteq X^2\setminus E.
$$
Define the auxiliary $R_k$ relations, $k\in\omega$, in
$2^{<\omega}$ by
\begin{align*}
s R_k t\iff & \lh(s)=\lh(t)\wedge (\forall i<k) s(i)=t(i)=0\\
&\wedge s(k)\neq t(k)\wedge (\forall i<\lh(s))(i>k\implies
s(i)=t(i)).
\end{align*}
We also let
$$
R=\bigcup_{k\in\omega} R_k.
$$
Fix winning $\Delta^1_1$ strategies for II in the strong Choquet
games on $(\baire,\tau)$ and $(\baire^2,\tau_2)$. Following
\cite[Lemma 6.3.10]{gao09}, it suffices to construct a scheme
consisting of $\tau$-basic open sets $(U_s)_{s\in 2^{<\omega}}$,
$(V_s)_{s\in 2^{<\omega}}$ that are subsets of $X$, points
$(x_s)_{s\in 2^{<\omega}}$ in $X$, and basic $\tau_2$-open sets
$(F_{s,t})_{sRt}$, $(E_{s,t})_{sRt}$ that are subsets of $X^2\cap
E$ such that
\begin{enumerate}[(i)]
\item $U_\emptyset=V_\emptyset=X$ \item For each $s\in
2^{<\omega}$ the following is a play according to II's winning
strategy:
$$
\begin{array}{rlllllllll}
\mathrm{I}  & x_{s\restrict 1}, U_{s\restrict 1} &                  & x_{s\restrict 2}, U_{s\restrict 2}  &                  & \cdots  &        & x_s,U_s\\
   &                                   &                  &                                     &                  &         &        &         &  \\
\mathrm{II} &                                   & V_{s\restrict 1}
& & V_{s\restrict 2} &         & \cdots &         & V_s
\end{array}
$$
\item $\diam(U_s)<2^{-\lh(s)}$ (with respect to the usual metric
on $\baire$). \item $U_{s\hat{\ }0}\times U_{s\hat{\ }1}\subseteq
C_{\lh(s)}$. \item If $\lh(s)=\lh(t)$ and $s R_k t$ then the
following is a play according to II's winning strategy in the
Choquet game on $(\baire^2,\tau_2)$:
$$
\begin{array}{rlllllllll}
\mathrm{I}  & (x_{s\restrict 1},x_{t\restrict 1}), F_{s\restrict 1,t\restrict k} &                   &                  & \cdots  &        & (x_s,x_t),F_{s,t}\\
   &                                   &                  &                                     &                  &         &        &         &  \\
\mathrm{II} &                                   & E_{s\restrict
1,t\restrict 1} &  &         & \cdots &         & E_{s,t}
\end{array}
$$
\item If $sRt$ then $\diam(F_{s,t})<2^{-\lh(s)}$ (with respect to
the usual metric on $\baire^2$).
\end{enumerate}
The construction of this scheme given in \cite{gao09} can easily
be carried out so that the function $s\mapsto (x_s, U_s, V_s)$ is
$\Delta^1_1$ relative to the set $C$. Thus $s\mapsto (x_s, U_s,
V_s)$ is $\Delta^1_1(O)$ and so the function defined by
$$
f(x)=y\iff y\in\bigcap_{n\in\omega} V_n
$$
is $\Delta^1_1(O)$. Finally, the arguments of \cite[p.
146--147]{gao09} show that $f$ is a reduction of $E_0$ to $E$.
\end{proof}

We will now show that Theorem \ref{silver} and \ref{hkl} are in
some sense optimal:

\begin{thm}
Let $z$ be a real in which $O$ is not hyperarithmetic. Then:
\begin{enumerate}[(i)]
\item There is a Hyp equivalence relation $E$ such that
$=_{\mathcal{P}(\omega)}\leq_{\Delta^1_1(O)} E$, but
$=_{\mathcal{P}(\omega)}\nleq_{\Delta^1_1(z)} E.$

\item There is a Hyp equivalence relation $E$ such that
$E_0\leq_{\Delta^1_1(O)} E$, but $E_0\nleq_{\Delta^1_1(z)} E$.
\end{enumerate}
\label{optimal}
\end{thm}
We need the following Lemma:

\begin{lemma}
Suppose $z$ is a real such that every non-empty $\Pi^0_1$ set in
$\baire$ contains a real hyperarithmetic in $z$. Then $O$ is
hyperarithmetic in z.
\end{lemma}
\begin{proof}
Let $C\subseteq\omega\times \baire$ be $\Pi^0_1$ and universal for
$\Pi^0_1$, and let
$$
\hat O=\{n: C_n=\emptyset\}.
$$
By our assumption,
$$
\omega \setminus \hat O=\{n:(\exists x\in\Delta_1^1(z)) x\in
C_n\},
$$
which is both $\Sigma^1_1$ and $\Pi^1_1(z)$. Thus $\hat O$ (and
therefore also Kleene's $O$) is $\Delta_1^1(z)$.
\end{proof}

\begin{proof}[Proof of Theorem \ref{optimal}]
Suppose $z$ is a real in which $O$ is not hyperarithmetic. Then by
the previous Lemma there is a non-empty $\Pi^0_1$ set
$F\subset\baire$ which does not contain any elements
hyperarithmetic in $z$, and in particular is uncountable. To prove
({\it i}), we let
$$
xEy\iff x,y\notin F\vee (x,y\in F\wedge x=y)
$$
Then $E$ has uncountably many classes. If $f:2^\omega\to\baire$
were a function witnessing that $=_{\mathcal{P}(\omega)}\leq_{\Delta^1_1(z)}
E$ then $f(\bar 0)\in F$ or $f(\bar 1)\in F$, which contradicts
that $F$ contains no real which is $\Delta^1_1(z)$.

To prove ({\it ii}), we instead define $E$ on $\baire\times
2^\omega$ by
$$
(x_0,y_0)E(x_1,y_1)\iff x_0,x_1\notin F\vee (x_0,x_1\in F\wedge
x_0=x_1\wedge y_0 E_0 y_1).
$$
Clearly $E$ is not smooth. If $E_0\leq_{\Delta^1_1(z)} E$ and
$f:2^\omega\to\baire\times 2^\omega$ witness this, then the
function
$$
\pi(x)=x_0\iff (\exists y_0) f(x)=(x_0,y_0)
$$
is constant on $E_0$ classes. It follows that there is some
$x_0\in F$ such that
$$
\{x\in 2^\omega: \pi(x)=x_0\},
$$
is comeagre. Since
$$
y\in\{x_0\}\iff \{x\in 2^\omega: \pi(x)=y\}\text{ is not meagre}
$$
it follows by \cite[Theorem 2.2.5]{kechris73} that $\{x_0\}$ is
$\Sigma^1_1(z)$, and so $x_0\in F$ is $\Delta^1_1(z)$, a
contradiction.
\end{proof}

\bigskip

\noindent \textbf{Remark.} Let $E$ be a Hyp equivalence relation
with $n\leq\omega$ many classes. Then by Proposition \ref{prop2},
$E\leq_{H} =_n$. On the other hand, since any $E$-class is Hyp it
must contain a real which is hyperarithmetic in $O$. Thus there is a
$\Delta_1^1(O)$ reduction of $=_n$ to $E$. We have the following
dichotomy:

\begin{thm}[The Finite Dichotomy Theorem] Let $E$ be a Hyp
equivalence relation on $\baire$. Then:
\begin{enumerate}[(a)]
\item If $n<\omega$ then either $E\leq_{H} =_n$ or
$=_{n+1}\leq_{\Delta^1_1(O)} E$.

\item Either there is $n<\omega$ such that $E\leq_{H}
=_n$, or $=_{\omega}\leq_{\Delta^1_1(O)} E$.
\end{enumerate}
\end{thm}
This Theorem is again optimal by an argument similar to that given
for Theorem \ref{optimal}.

\medskip

We conclude the paper with the following questions. The first
question seems related to \cite[Question 6.1.B]{jakelo02}.

\begin{quest}
 Can a Hyp equivalence relation $E$ be Borel reducible to $E_0$ but not Hyp reducible to $E_0$?
\end{quest}

\begin{quest}
 Are there any Hyp-degrees of Hyp equivalence relations other than $=_1$ that are comparable with all other Hyp-degrees?
\end{quest}

\begin{quest}
 What is the complexity of the first order theory of the partial order of Hyp-degrees of Hyp equivalence relations?
\end{quest}

\bibliographystyle{plain}
\bibliography{FFT}

\parbox{4.7in}{
{\sc
\noindent
Kurt G\"{o}del Research Center for Mathematical Logic \hfill \\
\hspace*{.1in}  University of Vienna \hfill \\
\hspace*{.2in}  W\"{a}hringer Stra{\ss}e 25 \hfill\\
\hspace*{.3in}  A-1090 Wien, Austria \hfill}\\
\hspace*{.045in} {\it E-mail: } \texttt{efokina\at {logic.univie.ac.at} }\hfill \\
\hspace*{.68in} \texttt{sdf\at {logic.univie.ac.at} }\hfill \\
\hspace*{.68in} \texttt{asger\at {logic.univie.ac.at}}\hfill

}
\end{document}